\documentclass[11pt,leqno]{amsart}

\usepackage{amscd}
\usepackage{amsmath}
\usepackage{amsfonts}
\usepackage{amssymb}
\usepackage{url}
\usepackage{graphicx}
\usepackage{pictexwd,dcpic}
\usepackage{dsfont}

\title{Galois realizability of groups of orders $p^5$ and $p^6$}
\author{Ivo M. Michailov}
\address{Faculty of Mathematics and Informatics, Shumen University "Episkop Konstantin Preslavski", Universitetska str. 115, 9700 Shumen, Bulgaria}
\email{ivo\_michailov@yahoo.com}
\date{\today}
\thanks{This work is partially supported by a project No RD-05-275/15.03.2012 of Shumen University.}
\begin{document}
\baselineskip 20pt
\begin{abstract}
Let $p$ be an odd prime, and let $k$ be an arbitrary field of
characteristic not $p$. In this article we determine the
obstructions for the realizability as Galois groups over $k$ of all
groups of orders $p^5$ and $p^6$, that have an abelian quotient
obtained by factoring out central subgroups of order $p$ or $p^2$.
These obstructions are decomposed as products of $p$-cyclic
algebras, provided that $k$ contains certain roots of unity.
\end{abstract}

\maketitle
\newcommand{\Gal}{{\rm Gal}}
\newcommand{\hw}{{\rm hw}}
\newcommand{\Ker}{{\rm Ker}}
\newcommand{\GL}{{\rm GL}}
\newcommand{\Br}{{\rm Br}}
\newcommand{\lcm}{{\rm lcm}}
\newcommand{\ord}{{\rm ord}}
\newcommand{\res}{{\rm res}}
\newcommand{\cor}{{\rm cor}}
\newcommand{\ch}{{\rm char}}
\newcommand{\tr}{{\rm tr}}
\newcommand{\Hom}{{\rm Hom}}
\newcommand{\ind}{{\rm ind}}
\newcommand{\Pin}{{\rm Pin}}
\newcommand{\Spin}{{\rm Spin}}
\newcommand{\attop}[2]{\genfrac{}{}{0pt}{1}{#1}{#2}}
\renewcommand{\thefootnote}{\fnsymbol{footnote}}
\numberwithin{equation}{section}

\makeatletter

%--------------------------------------Section 1-------------------------------------------------
\section{Introduction}
\label{1}

Let $p$ be an odd prime, and let $k$ be an arbitrary field of
characteristic not $p$. In this article we investigate the
realizability of a number of groups of orders $p^5$ and $p^6$ as
Galois groups over $k$, provided that $k$ contains certain roots of
unity. Recall that given a group $G$, the inverse problem of Galois
theory asks whether or not there exists a Galois field extension $K$
of $k$ such that the Galois group $\Gal(K/k)$ is isomorphic to $G$.

In recent years, extensive research has been done relating to the
description of necessary and sufficient conditions for the
realizability of small $p$-groups, especially for $p=2$. These
conditions are often expressed as the splitting of certain elements
in the Brauer group $\Br(k)$, called \emph{obstructions}. The most
important part in the investigations is decomposing the obstructions
as cyclic algebras (or, quaternion algebras for $p=2$).

Ledet \cite{Le-95} found the decompositions of the obstructions as
quaternion algebras of all non abelian groups of order $2^n$ for
$n\leq 4$. The obstructions for the realizability of groups of
orders $32$ and $64$ are calculated in a number of papers e.g.
\cite{Mi-32,GS-D4,GS-32,Mi-coh,GS-64}. Massy \cite{Ma-87}
investigated the realizability of the two non abelian groups of
order $p^3$. Michailov \cite{Mi-p4} determined the obstructions to
the realizability of four non abelian groups of order $p^4$. An
extensive survey of the realizability of $p$-groups has been done
recently by Michailov and Ziapkov \cite{MZ-pGr}.

The present paper is focused on the obstructions for the
realizability of all groups of orders $p^5$ and $p^6$ that have an
abelian quotient obtained by factoring out central subgroups of
order $p$ or $p^2$. In order to achieve this goal, we develop
several theoretic criteria given in Section \ref{2}. In Section
\ref{3} we list the groups that fulfill our condition. We use the
classification of the groups of orders $p^5$ and $p^6$ made by James
\cite{Ja}. In Section \ref{4} we determine the obstructions for the
realizability of $18+(p-1)$ groups of order $p^5$ and more than
$71+4(p-1)$ groups of order $p^6$, and write them in six tables.

\bigskip

%--------------------------------------Section 2-------------------------------------------------
\section{Cohomological criteria}
\label{2}

Let $k$ be an arbitrary field and let $G$ be a non simple group.
Assume that $A$ is a normal subgroup of $G$. Then the realizability
of the quotient group $F=G/A$ as a Galois group over $k$ is a
necessary condition for the realizability of $G$ over $k$. In this
way arises the next generalization of the inverse problem in Galois
theory -- the embedding problem of fields.

Let $K/k$ be a Galois extension with Galois group $F$, and let
\begin{equation}\label{e2.1}
1 \longrightarrow A \longrightarrow
G\overset{\alpha}{\longrightarrow} F \longrightarrow 1,
\end{equation}
be a group extension, i.e., a short exact sequence. Solving {\it the
embedding problem} related to $K/k$ and \eqref{e2.1} consists of
determining whether or not there exists a Galois algebra (called
also a {\it weak} solution) or a Galois extension (called a {\it
proper} solution) $L$, such that $K$ is contained in $L$, $G$ is
isomorphic to $\Gal(L/k)$, and the homomorphism of restriction to
$K$ of the automorphisms from $G$ coincides with $\alpha$. We denote
the so formulated embedding problem by $(K/k, G, A)$. We call the
group $A$ the {\it kernel} of the embedding problem.

Now, let $k$ be an arbitrary field of characteristic not $p$,
containing a primitive $p^n$th root of unity $\zeta_{p^n}$ for
$n\in\mathds N$, and put $\mu_{p^n}=\langle\zeta_{p^n}\rangle$. Let
$K$ be a Galois extension of $k$ with Galois group $F$. Consider a
non split group extension
\begin{equation}\label{e3.1}
1\longrightarrow \langle\varepsilon\rangle\longrightarrow
G\longrightarrow F\longrightarrow 1,
\end{equation}
where $\varepsilon$ is a central element of order $p^n$ in $G$. We
are going to identify the groups $\langle\varepsilon\rangle$ and
$\mu_{p^n}$, since they are isomorphic as $F$-modules.

Assume that $c\in H^2(F,\mu_{p^n})$ is the $2$-coclass corresponding
to the group extension \eqref{e3.1} and denote by $\Omega_k$ the
Galois group of the algebraic separable closure $\bar k$ over $k$.
\emph{The obstruction} to the embedding problem $(K/k,G,\mu_{p^n})$
we call the image of $c$ under the inflation map
$\inf_F^{\Omega_k}:H^2(F,\mu_{p^n})\to H^2(\Omega_k,\mu_{p^n})$.

Note that we have the standard isomorphism of
$H^2(\Omega_k,\mu_{p^n})$ with the ${p^n}$-torsion in the Brauer
group of $k$ induced by applying $H^\ast(\Omega_k,\cdot)$ to the
$p^n$-th power exact sequence of $\Omega_k$-modules
$1\longrightarrow\mu_{p^n}\longrightarrow\bar
k^\times\longrightarrow\bar k^\times\longrightarrow 1$. In this way,
the obstruction equals the equivalence class of the crossed product
algebra $(F,K/k,\bar c)$ for any $\bar c\in c$. Hence we may
identify the obstruction with a Brauer class in $\Br_{p^n}(k)$.

Note that we have an injection $\mu_{p^n}\hookrightarrow K^\times$,
which induces a homomorphism $\nu:H^2(F,\mu_{p^n})\to
H^2(F,K^\times)$. Then the obstruction is equal to $\nu(c)$, since
there is an isomorphism between the relative Brauer group ${\rm
Br}(K/k)$ and the group $H^2(F,K^\times)$.

More generally, the following result holds.

\newtheorem{t3.1}{Theorem}[section]
\begin{t3.1}\label{t3.1}
{\rm (\cite{Ki,MZ},)} Let $n\geq 1$, and let $c$ be the $2$-coclass
in $H^2(F,\mu_{p^n})$, corresponding to the non split central group
extension \eqref{e3.1}. Then the embedding problem
$(K/k,G,\mu_{p^n})$ is weakly solvable if and only if $\nu(c)=1$. If
$n=1$ or $\mu_{p^n}$ is contained in the Frattini subgroup $\Phi(G)$
of $G$ (for $n>1$), then the condition $\nu(c)=1$ is sufficient also
for the proper solvability of the problem $(K/k,G,\mu_{p^n})$ (see
\cite[\S 1.6, Cor. 5]{ILF}).
\end{t3.1}

Henceforth, embedding problems of the kind $(K/k,G,\mu_{p^n})$ we
will call $\mu_{p^n}$-embedding problems. We are going to consider
first the case $n=1$.

From the well-known Merkurjev-Suslin Theorem \cite{MeS} it follows
that the obstruction to any  $\mu_p$-embedding problem is equal to a
product of classes of $p$-cyclic algebras. The explicit computation
of these $p$-cyclic algebras, however, is not a trivial task. We are
going to discuss the methods for achieving this goal.

We denote by $(a, b;\zeta)$ the equivalence class of the $p$-cyclic
algebra which is generated by $i_1$ and $i_2$, such that $i_1^p=b,
i_2^p=a$ and $i_1i_2=\zeta i_2i_1$. For $p=2$ we have the quaternion
class $(a,b;-1)$, commonly denoted by $(a,b)$.

In 1987 Massy \cite{Ma-87} obtained a formula for the decomposition
of the obstruction in the case when $F=\Gal(K/k)$ is isomorphic to
$(C_p)^n$, the elementary abelian $p$-group.

\newtheorem{t3.2}[t3.1]{Theorem}
\begin{t3.2}\label{t3.2}
{\rm (\cite[Th\'eor\`eme 2]{Ma-87},\cite[Cor. 6.1.6]{Le-05})} Let
$K/k=k(\root p\of {a_1},\root p\of {a_2},\dots,\root p\of {a_n})/k$
be a $(C_p)^n$ extension, and let
$\sigma_1,\sigma_2,\dots,\sigma_n\in {\rm Gal}(K/k)$ be given by
$\sigma_i(\root p\of {a_j})/\root p\of {a_j}=\zeta^{\delta_{ij}}$
($\delta_{ij}$ is the Kronecker delta). Let
\begin{equation*}
1\longrightarrow \mu_p\longrightarrow G\longrightarrow {\rm
Gal}(K/k)\longrightarrow 1
\end{equation*}
be a non split central extension, and choose pre-images
$s_1,s_2,\dots,s_n\in G$ of $\sigma_1,\sigma_2,\dots,\sigma_n$.
Define $d_i (1\leq i\leq n)$ by $s_i^p=\zeta^{d_i}$, and $d_{ij}
(i<j)$  by $s_is_j=\zeta^{d_{ij}}s_js_i$. Then the obstruction to
the proper solvability of the embedding problem $(K/k,G,\mu_p)$ is
\begin{equation*}
\prod_{i=1}^{n}(a_i,\zeta;\zeta)^{d_i}\prod_{i<j}(a_j,a_i;\zeta)^{d_{ij}}.
\end{equation*}
\end{t3.2}

Michailov \cite{Mi-mod,Mi-p4} obtained a formula for the
decomposition of the obstruction in the case when the quotient $F$
has a direct factor $C_p$.

Let $H$ be a $p$-group and let
\begin{equation}\label{e3.2}
1\longrightarrow \mu_p\longrightarrow
G\overset{\pi}{\longrightarrow} F\cong H\times C_p\longrightarrow 1
\end{equation}
be a non split central group extension with characteristic
$2$-coclass $\gamma\in H^2(H\times C_p,C_p)$. By $\res_H\gamma$ we
denote the $2$-coclass of the group extension
\begin{equation*}
1\longrightarrow \mu_p\longrightarrow
\pi^{-1}(H)\overset{\pi}{\longrightarrow} H\longrightarrow 1.
\end{equation*}
Let $\sigma_1,\sigma_2,\dots,\sigma_m$ be a minimal generating set
for the maximal elementary abelian quotient group of $H$; and let
$\tau$ be the generator of the direct factor $C_p$. Finally, let
$s_1,s_2,\dots,s_m,t\in G$ be the pre-images of
$\sigma_1,\sigma_2,\dots,\sigma_m,\tau$, such that $t^p=\zeta^j$ and
$ts_i=\zeta^{d_i}s_it$, where $i\in\{1,2,\dots,m\};
j,d_i\in\{0,1,\dots,p-1\}$.

\newtheorem{t3.3}[t3.1]{Theorem}
\begin{t3.3}\label{t3.3}
{\rm(\cite[Theorem 4.1]{Mi-mod},\cite[Theorem 2.1]{Mi-p4})} Let
$K/k$ be a Galois extension with Galois group $H$ and let
$L/k=K(\root p\of b)/k$ be a Galois extension with Galois group
$H\times C_p$ ($b\in k^\times\setminus k^{\times p}$). Choose
$a_1,a_2,\dots,a_m\in k^\times$, such that $\sigma_k\root p\of
{a_i}=\zeta^{\delta_{ik}}\root p\of {a_i}$ ($\delta_{ik}$ is the
Kronecker delta). Then the obstruction to the proper solvability of
the embedding problem $(L/k,G,\mu_p)$ is
\begin{equation*}
[K,H,\res_H\gamma]\left(b,\zeta^j\prod_{i=1}^{m}a_i^{d_i};\zeta\right).
\end{equation*}
\end{t3.3}

Ledet describes in his book \cite{Le-05} a more general formula for
the decomposition of the obstruction of $\mu_p$-embedding problems
with finite group $F$ isomorphic to a direct product of two groups.

Let $G$ be arbitrary finite group, and let $p$ be a prime divisor of
$\ord(G)$. Define $\mathcal O^p(G)$ as the subgroup of $G$ generated
by all elements of order prime to $p$. It is clear that $\mathcal
O^p(G)$ is the intersection of all normal subgroups in $G$ of
$p$-power index.

\newtheorem{t3.4}[t3.1]{Theorem}
\begin{t3.4}\label{t3.4}
{\rm(\cite[Theorem 6.1.4]{Le-05}} Let $L/k$ be a $N\times H$
extension, where $N$ and $H$ are finite groups. Let
\begin{equation*}
1\longrightarrow \mu_p\longrightarrow G\longrightarrow N\times
H\longrightarrow 1
\end{equation*}
be a non split central group extension with cohomology class
$\gamma\in H^2(N\times H,\mu_p)$. Let $K'/k$ and $K/k$ be the
subextensions corresponding to the factors $N$ and $H$. (I.e.,
$K'=L^H,K=L^N$.) Let $\sigma_1,\sigma_2,\dots,\sigma_m$ and
$\tau_1,\tau_2,\dots,\tau_n$ represent minimal generating sets for
the groups $N/\mathcal O^p(N)$ and $H/\mathcal O^p(H)$, and choose
$a_1,a_2,\dots,a_m,b_1,\dots b_n\in k^\times$, such that $\root p\of
{a_i}\in K^\times,\root p\of {b_i}\in K'^\times, \sigma_\kappa\root
p\of {a_i}=\zeta^{\delta_{i\kappa}}\root p\of {a_i}$ and
$\tau_\ell\root p\of {b_i}=\zeta^{\delta_{i\ell}}\root p\of {b_i}$
($\delta$ is the Kronecker delta). Finally, let
$s_1,\dots,s_m;t_1,\dots,t_n\in G$ be the pre-images of
$\sigma_1,\dots,\sigma_m;\tau_1,\dots,\tau_n$, and let
$d_{ij}\in\{0,\dots,p-1\}$ be given by
$t_js_i=\zeta^{d_{ij}}s_it_j$.

Then the obstruction to the proper solvability of the embedding
problem $(L/k,G,\mu_p)$ given by $\gamma$ is
\begin{equation*}
[K,N,\res_N\gamma]\cdot [K',H,\res_H\gamma]\cdot
\prod_{i,j}(b_j,a_i;\zeta)^{d_{ij}}.
\end{equation*}
\end{t3.4}

Now, we will prove our first main result that will allow us to
decompose the obstruction when the quotient is an arbitrary abelian
$p$-group.

\newtheorem{c3.6}[t3.1]{Theorem}
\begin{c3.6}\label{c3.6}
Let $L/k$ be an $H\cong\prod_{i=1}^t C_{p^{n_i}}$ extension for some
natural numbers $n_1\leq n_2\leq\cdots\leq n_t$. Let
\begin{equation*}
1\longrightarrow \mu_p\longrightarrow G\longrightarrow H\cong
\prod_{i=1}^t C_{p^{n_i}} \longrightarrow 1
\end{equation*}
be a non split central group extension with cohomology class
$\gamma\in H^2(H,\mu_p)$. Let $K_i/k$ be the subextension
corresponding to the factor $C_{p^{n_i}}$ for $i=1,\dots,t$. (I.e.,
$K_i$ is the fixed subfield of $\prod_{j\ne i}C_{p^{n_j}}$.) Let
$\sigma_i$ be the generator of $C_{p^{n_i}}$ for $i=1,\dots,t$, and
choose $a_i\in k^\times$, such that $\root p\of {a_i}\in K_i^\times$
and $\sigma_j\root p\of {a_i}=\zeta^{\delta_{ij}}\root p\of {a_i}$
($\delta$ is the Kronecker delta). Let $s_1,\dots,s_t$ be the
pre-images of $\sigma_1,\dots,\sigma_t$, let
$d_{ij}\in\{0,\dots,p-1\}$ be given by $s_i s_j=\zeta^{d_{ji}}s_j
s_i$, and let $s_i^{p^{n_i}}=\zeta^{m_i}$ for $i=1,\dots,t;
m_i\in\{0,\dots,p-1\}$.

Finally, define $r=\max\{i:m_i>0\},n=n_r,A=\{i:n_i=n,m_i>0\}$, and
assume that $k$ contains $\zeta_{p^n}$, a primitive $p^n$-th root of
unity. Then the obstruction to the proper solvability of the
embedding problem $(L/k,G,\mu_p)$ given by $\gamma$ is
\begin{equation*}
\prod_{i\in A} (a_i,\zeta_{p^n}^{m_i};\zeta)\cdot
\prod_{i<j}(a_j,a_i;\zeta)^{d_{ij}}.
\end{equation*}
\end{c3.6}
\begin{proof}
From Theorem \ref{t3.4} by induction we obtain that the obstruction
to the proper solvability of the embedding problem $(L/k,G,\mu_p)$
given by $\gamma$ is
\begin{equation*}
\prod_{i=1}^{t} [K_i,C_{p^{n_i}},\res_{C_{p^{n_i}}}\gamma]\cdot
\prod_{i<j}(a_j,a_i;\zeta)^{d_{ij}}.
\end{equation*}
Consider the cyclic algebra $[K_r,\sigma_r,\zeta^{m_r}]$. It is
generated by the field $K_r$ such that $[K_r:k]=p^n$ and an element
$s_r$ such that $s_r^{p^n}=\zeta^{m_r}$. Denote by $k_r/k=k(\root
p\of {a_r})/k$ the subextension that is contained in $K_r/k$, and
consider the cyclic subalgebra
$[k_r,\sigma_r\vert_{k_r},\zeta_{p^n}^{m_r}]$ of
$[K_r,\sigma_r,\zeta^{m_r}]$, where
$\sigma_r^p\vert_{k_r}=\zeta_{p^n}^{m_r}$. Then, according to
\cite[Corollary 15.1b]{Pi}, we have
$[K_r,\sigma_r,\zeta^{m_r}]=[K_r,\sigma_r,\zeta_{p^n}^{m_rp^n/p}]=[k_r,\sigma_r\vert_{k_r},\zeta_{p^n}^{m_r}]=(a_r,\zeta_{p^n}^{m_r};\zeta)\in
\Br(k)$. Note that this fact is a variation of the well known
Albert's Theorem \cite[Theorem 11, p. 207]{Al}, namely that
$[K_r,\sigma_r,\zeta^{m_r}]$ is split if and only if
$\zeta_{p^n}^{m_r}\in N_{k_r/k}(k_r^\times)$, where by $N_{k_r/k}$
we denote the norm map.

Now, let $i$ be an arbitrary integer such that $n_i<n,m_i>0$, and
put $m=n_i$. As we have just shown,
$[K_i,\sigma_i,\zeta^{m_i}]=(a_i,\zeta_{p^m}^{m_i};\zeta)$. Since
$\zeta_{p^m}=\zeta_{p^n}^{p^{n-m}}$ for $n-m\geq 1$, we obtain that
$(a_i,\zeta_{p^m}^{m_i};\zeta)=1\in\Br(k)$. In this way, we obtain
the formula given in the statement.
\end{proof}

With the aid of the latter result we will determine in Section
\ref{4} the obstructions to the realizability of all $9$ groups of
order $p^5$ and all $15$ groups of order $p^6$ that have an abelian
quotient group, and that are not a direct product of smaller groups.

We extend this method to one applying to an additional $9+(p-1)$
groups of order $p^5$ and more than $53+4(p-1)$ \footnote{The actual
number of the groups $\Phi_{15}(2211)b_{r,s}$ seems difficult to
calculate, since the values of $s$ depend in a complicated way on
the values of $r$. That is why we counted only the values of $r$.}
groups of order $p^6$ with the following property: there exist two
disjoint central subgroups $N_1$ and $N_2$ of order $p$, such that
the quotient group obtained by factoring out $N_1N_2$ is abelian.
Any group that has two disjoint central subgroups $N_1$ and $N_2$ is
a pullback. Pullbacks of orders $16,32$ and $64$ are considered in
\cite{Le-95,Mi-32,GS-64}.

Namely, let $\varphi' : G'\rightarrow F$ and $\varphi'' :
G''\rightarrow F$ be homomorphisms with kernels $N'$ and,
respectively, $N''$. {\it The pullback} of the pair of homomorphisms
$\varphi'$ and $\varphi''$ is the subgroup in $G'\times G''$ of all
pairs $(\sigma',\sigma'')$, such that
$\varphi'(\sigma')=\varphi''(\sigma'')$. The pullback is denoted by
$G'\curlywedge G''$. It is also called the direct product of the
groups $G'$ and $G''$ with amalgamated quotient group $F$ and
denoted by $G'*_F G''$.

Now, let $N_1=N'\times \{1\}$ and $N_2=\{1\}\times N''$. Then $N_1$
and $N_2$ are normal subgroups of $G'\curlywedge G''$, such that
$N_1\cap N_2=\{1\}$. The converse is also true (see \cite{ILF}, I,
\S 12):

\newtheorem{l4.1}[t3.1]{Lemma}
\begin{l4.1}\label{l4.1}
Let $N_1$ and $N_2$ be two normal subgroups of the group $G$, such
that $N_1\cap N_2=\{1\}$. Then $G$ is isomorphic to the pullback
$(G/N_1)\curlywedge (G/N_2)$.
\end{l4.1}

The application to embedding problems is given by:

\newtheorem{t4.0}[t3.1]{Theorem}
\begin{t4.0}\label{t4.0}
{\rm (\cite[Theorem 1.12]{ILF})} Let $K/k$ be a Galois extension
with Galois group $F$. In the notations of the lemma, let $F\cong
G/N_1N_2$ and $G\cong (G/N_1)\curlywedge (G/N_2)$. Then the
embedding problem $(K/k, G, N_1N_2)$ is solvable if and only if the
embedding problems $(K/k, G/N_1, N_2)$ and $(K/k, G/N_2, N_1)$ are
solvable.
\end{t4.0}

Next, we are going to prove a result for certain
$\mu_{p^n}$-embedding problems for an arbitrary $n\in \mathds N$. We
will find a decomposition of the obstruction as a product of
$p^n$-cyclic algebras. These algebras are similar to the $p$-cyclic
algebras, mentioned earlier in this Section. Namely, we denote by
$(a, b;\zeta_{p^n})$ the equivalence class of the cyclic algebra
which is generated by $i_1$ and $i_2$, such that $i_1^{p^n}=b,
i_2^{p^n}=a$ and $i_1i_2=\zeta_{p^n} i_2i_1$. Of course, we assume
again that $k$ contains $\zeta_{p^n}$, a primitive $p^n$-th root of
unity. For more details about these algebras we refer the reader to
\cite[\S 15]{Pi}.

\newtheorem{t3.7}[t3.1]{Theorem}
\begin{t3.7}\label{t3.7}
Let $k$ contain a primitive $p^n$-th root of unity $\zeta_{p^n}$,
and let $L/k=k(\root{p^n}\of {a_1},\dots,$ $\root{p^n}\of {a_m})/k$
be an arbitrary $(C_{p^n})^m$ extension for some $m,n\in\mathds N$.
Let
\begin{equation*}
1\longrightarrow \mu_{p^n}\longrightarrow G\longrightarrow
(C_{p^n})^m\longrightarrow 1
\end{equation*}
be a non split central group extension, let
$\sigma_1,\sigma_2,\dots,\sigma_m$ be the generators of
$(C_{p^n})^m$, and let $s_1,s_2,\dots,s_m\in G$ be their pre-images
such that $s_i^{p^n}=\zeta_{p^n}^{j_i}$ and
$s_js_i=\zeta_{p^n}^{d_{ij}}s_is_j$, where $i\in\{1,2,\dots,m\};
j_i,d_{ij}\in\{0,1,\dots,p^n-1\}$ and $i<j$. Assume that
$\sigma_j\root{p^n}\of {a_i}=\zeta_{p^n}^{\delta_{ij}}\root{p^n}\of
{a_i}$ ($\delta_{ij}$ is the Kronecker delta). Then the obstruction
to the weak solvability of  the embedding problem
$(L/k,G,\mu_{p^n})$ is
\begin{equation*}
\prod_{i=1}^m (a_i,\zeta_{p^n}^{j_i};\zeta_{p^n})\cdot
\prod_{i<j}(a_j,a_i;\zeta_{p^n})^{d_{ij}}.
\end{equation*}
\end{t3.7}
\begin{proof}
Let $\mathcal A=(L,(C_{p^n})^m,\zeta_{p^n})$ be the crossed product
algebra related to the embedding problem $(L/k,G,\mu_{p^n})$. Denote
$H=\langle\sigma_1,\dots,\sigma_{m-1}\rangle\cong (C_{p^n})^{m-1}$
and $K/k=k(\root{p^n}\of {a_1},\dots,$ $\root{p^n}\of {a_{m-1}})/k$.
The crossed product algebra $\mathcal B=(K,H,\zeta_{p^n})$ is
included in $\mathcal A$, therefore $\mathcal A$ is a tensor product
of $B$ and the centralizer of $\mathcal B$ in $\mathcal A: \mathcal
A=\mathcal B\otimes_kC_A(B)$. Now, consider the subalgebra $\mathcal
C=k[\root{p^n}\of{a_m},\left(\prod_{i=1}^{m-1}\root{p^n}\of
{a_i}^{d_{im}}\right)s_m]$ in $A$. Since
$s_m^{p^n}=\zeta_{p^n}^{j_m},s_m\root{p^n}\of{a_m}=\zeta_{p^n}\root{p^n}\of{a_m}s_m$,
and $s_m$ commutes with $\prod_{i=1}^{m-1}\root{p^n}\of
{a_i}^{d_{im}}$, we have the isomorphism $\mathcal C\cong
(a_m,\zeta_{p^n}^{j_m}\prod_{i=1}^{m-1} a_i^{d_{im}};\zeta_{p^n})$.

Next, we will show that $\mathcal C$ is in fact the centralizer
$C_{\mathcal A}(\mathcal B)$. Indeed, for $1\leq\kappa\leq m-1$ we
have
\begin{equation*}\begin{split}
&s_\kappa\left(\prod_{i=1}^{m-1}\root{p^n}\of
{a_i}^{d_{im}}\right)s_m=\left(\prod_{i=1}^{m-1}\zeta_{p^n}^{\delta_{i\kappa}d_{im}}\root{p^n}\of
{a_i}^{d_{im}}\right)s_\kappa
s_m=\left(\prod_{i=1}^{m-1}\zeta_{p^n}^{\delta_{i\kappa}d_{im}}\root{p^n}\of
{a_i}^{d_{im}}\right)\zeta_{p^n}^{-d_{\kappa m}}s_m s_\kappa\\
&=\zeta_{p^n}^{-d_{\kappa
m}}\prod_{i=1}^{m-1}\zeta_{p^n}^{\delta_{i\kappa}d_{im}}\left(\prod_{i=1}^{m-1}\root{p^n}\of
{a_i}^{d_{im}}\right)s_m
s_\kappa=\left(\prod_{i=1}^{m-1}\root{p^n}\of
{a_i}^{d_{im}}\right)s_m s_\kappa,
\end{split}\end{equation*}
since $\sum_{i=1}^{m}\delta_{i\kappa}d_{im}=d_{\kappa m}$.
Therefore, $[\mathcal A]=[\mathcal
B](a_m,\zeta_{p^n}^{j_m}\prod_{i=1}^{m-1}
a_i^{d_{im}};\zeta_{p^n})$, and the theorem follows by induction.
\end{proof}

With the aid of the latter Theorem we will find in Section \ref{4}
the obstructions to 3 groups of order $p^6$ that can not be treated
with the previous criteria.

%--------------------------------------Section 3-------------------------------------------------
\section{The groups}
\label{3}

In this Section we give a list of groups of orders $p^5$ and $p^6$
that will be investigated for realizability in Section \ref{4}. We
use the classification made by R. James \cite{Ja}. The groups in
James' list are collected in a number of so called \emph{isoclinism}
families. Two groups $G,H$ with centers $Z(G),Z(H)$ and derived
groups $G',H'$ are said to be \emph{isoclinic} if there exist
isomorphisms $\theta:G/Z(G)\to H/Z(H)$ and $\phi: G'\to H'$ such
that $\phi([\alpha,\beta])=[\alpha',\beta']$ for all
$\alpha,\beta\in G$, where $\alpha' Z(H)=\theta(\alpha Z(G))$ and
$\beta'(Z(H))=\theta(\beta Z(G))$. This relation is well defined and
is in fact an equivalence relation. The equivalence classes are
called \emph{isoclinism families}. A family of $p$-groups will be
denoted by $\Phi_s$ where $p$ is an arbitrary prime and $s$ is an
integer.

\medskip
{\bf Standing notations.}
$[\alpha,\beta]=\alpha^{-1}\beta^{-1}\alpha\beta$. For economy of
space, all relations of the form $[\alpha,\beta]=1$ (with
$\alpha,\beta$ generators) have been omitted from the list and
should be assumed when reading the list. Throughout, $\nu$ denotes
the smallest positive integer which is a non-quadratic residue$\pmod
p$ and $g$ denotes the smallest positive integer which is a
primitive root$\pmod p$.
\medskip

By studying the list given in \cite{Ja}, we conclude that the groups
of order $p^5$ having an abelian quotient (obtained by factoring out
a central cyclic group of order $p$) are precisely those from
families $(2)$ and $(5)$ listed below. Of course, we will omit from
our lists the groups of the kind $G\times H$, since the
realizability of these groups depends only on the realizability of
the direct factors $G$ and $H$.

\begin{enumerate}
\item[{\rm (2)}]
{\allowdisplaybreaks
\begin{align*}
\Phi_2(41)=& \langle \alpha,\alpha_1,\alpha_2:[\alpha_1,\alpha]=\alpha^{p^3}=\alpha_2,~\alpha_1^p=\alpha_2^p=1\rangle, \\
\Phi_2(32)a1=& \langle \alpha,\alpha_1,\alpha_2:[\alpha_1,\alpha]=\alpha^{p^2}=\alpha_2,\alpha_1^{p^2}=\alpha_2^p=1\rangle, \\
\Phi_2(32)a2=& \langle \alpha,\alpha_1,\alpha_2:[\alpha_1,\alpha]=\alpha_1^p=\alpha_2,~\alpha^{p^3}=\alpha_2^p=1\rangle, \\
\Phi_2(311)b=& \langle \alpha,\alpha_1,\alpha_2,\gamma:[\alpha_1,\alpha]=\gamma^{p^2}=\alpha_2,~\alpha^p=\alpha_1^p=\alpha_2^p=1\rangle, \\
\Phi_2(311)c=& \langle \alpha,\alpha_1,\alpha_2:[\alpha_1,\alpha]=\alpha_2,~\alpha^{p^3}=\alpha_1^p=\alpha_2^p=1\rangle, \\
\Phi_2(221)c=& \langle \alpha,\alpha_1,\alpha_2,\gamma:[\alpha_1,\alpha]=\gamma^p=\alpha_2,~\alpha^{p^2}=\alpha_1^p=\alpha_2^p=1\rangle, \\
\Phi_2(221)d=& \langle
\alpha,\alpha_1,\alpha_2:[\alpha_1,\alpha]=\alpha_2,~\alpha^{p^2}=\alpha_1^{p^2}=\alpha_2^p=1\rangle;
\end{align*}}
\item[{\rm (5)}]
{\allowdisplaybreaks
\begin{align*}
\Phi_5(2111)=& \langle \alpha_1,\alpha_2,\alpha_3,\alpha_4,\beta:[\alpha_1,\alpha_2]=[\alpha_3,\alpha_4]=\alpha_1^p=\beta,~\alpha_2^p=\alpha_3^p=\alpha_4^p=\beta^p=1\rangle, \\
\Phi_5(1^5)=& \langle
\alpha_1,\alpha_2,\alpha_3,\alpha_4,\beta:[\alpha_1,\alpha_2]=[\alpha_3,\alpha_4]=\beta,~\alpha_1^p=\alpha_2^p=\alpha_3^p=\alpha_4^p=\beta^p=1\rangle.
\end{align*}}
\end{enumerate}
\medskip

Next, one can see that the pullbacks of order $p^5$ with an abelian
quotient\footnote{I.e., having two disjoint central subgroups $N_1$
and $N_2$ of order $p$, such that the quotient group obtained by
factoring out $N_1N_2$ is abelian.} are precisely those from family
$(4)$.

\begin{enumerate}
\item[{\rm (4)}]
{\allowdisplaybreaks
\begin{align*}
\Phi_4(221)a=& \langle \alpha,\alpha_1,\alpha_2,\beta_1,\beta_2:[\alpha_i,\alpha]=\beta_i,~\alpha^p=\beta_2,~\alpha_1^p=\beta_1,\alpha_2^p=\beta_i^p=1~ (i=1,2)\rangle, \\
\Phi_4(221)b=& \langle \alpha,\alpha_1,\alpha_2,\beta_1,\beta_2:[\alpha_i,\alpha]=\beta_i,~\alpha^p=\beta_2,~\alpha_2^p=\beta_1,\alpha_1^p=\beta_i^p=1~ (i=1,2)\rangle, \\
\Phi_4(221)c=& \langle \alpha,\alpha_1,\alpha_2,\beta_1,\beta_2:[\alpha_i,\alpha]=\beta_i=\alpha_i^p,~\alpha^p=\beta_i^p=1~ (i=1,2)\rangle, \\
\Phi_4(221)d_r=& \langle \alpha,\alpha_1,\alpha_2,\beta_1,\beta_2:[\alpha_i,\alpha]=\beta_i,~\alpha_1^p=\beta_1^\kappa,~\alpha_2^p=\beta_2,~\alpha^p=\beta_i^p=1~ (i=1,2)\rangle, \\
&\text{where}\ \kappa=g^r\ \text{for}\ r=1,2,\dots,\frac{1}{2}(p-1),\\
\Phi_4(221)e=& \langle \alpha,\alpha_1,\alpha_2,\beta_1,\beta_2:[\alpha_i,\alpha]=\beta_i,~\alpha_1^p=\beta_2^{-1/4},~\alpha_2^p=\beta_1\beta_2,\alpha^p=\beta_i^p=1~ (i=1,2)\rangle, \\
\Phi_4(221)f_0=& \langle \alpha,\alpha_1,\alpha_2,\beta_1,\beta_2:[\alpha_i,\alpha]=\beta_i,~\alpha_1^p=\beta_2,~\alpha_2^p=\beta_1^\nu,\alpha^p=\beta_i^p=1~ (i=1,2)\rangle, \\
\Phi_4(221)f_r=& \langle \alpha,\alpha_1,\alpha_2,\beta_1,\beta_2:[\alpha_i,\alpha]=\beta_i,~\alpha_1^p=\beta_2^\kappa,~\alpha_2^p=\beta_1\beta_2,\alpha^p=\beta_i^p=1~ (i=1,2)\rangle, \\
&\text{where}\ 4\kappa=g^{2r+1}\ \text{for}\ r=1,2,\dots,\frac{1}{2}(p-1),\\
\Phi_4(2111)a=& \langle \alpha,\alpha_1,\alpha_2,\beta_1,\beta_2:[\alpha_i,\alpha]=\beta_i,~\alpha^p=\beta_2,~\alpha_i^p=\beta_i^p=1~ (i=1,2)\rangle, \\
\Phi_4(2111)b=& \langle \alpha,\alpha_1,\alpha_2,\beta_1,\beta_2:[\alpha_i,\alpha]=\beta_i,~\alpha_1^p=\beta_1,~\alpha^p=\alpha_2^p=\beta_i^p=1~ (i=1,2)\rangle, \\
\Phi_4(2111)c=& \langle \alpha,\alpha_1,\alpha_2,\beta_1,\beta_2:[\alpha_i,\alpha]=\beta_i,~\alpha_2^p=\beta_1,~\alpha^p=\alpha_1^p=\beta_i^p=1~ (i=1,2)\rangle, \\
\Phi_4(1^5)=& \langle
\alpha,\alpha_1,\alpha_2,\beta_1,\beta_2:[\alpha_i,\alpha]=\beta_i,~\alpha^p=\alpha_i^p=\beta_i^p=1~
(i=1,2)\rangle.
\end{align*}}
\end{enumerate}

Similarly, the groups of order $p^6$ having an abelian quotient
(obtained by factoring out a central cyclic group of order $p$) are
precisely those from families $(2)$ and $(5)$.

\begin{enumerate}
\item[{\rm (2)}]
{\allowdisplaybreaks
\begin{align*}
\Phi_2(51)=& \langle \alpha,\alpha_1,\alpha_2:[\alpha_1,\alpha]=\alpha^{p^4}=\alpha_2,~\alpha_1^p=\alpha_2^p=1\rangle, \\
\Phi_2(42)a1=& \langle \alpha,\alpha_1,\alpha_2:[\alpha_1,\alpha]=\alpha^{p^3}=\alpha_2,\alpha_1^{p^2}=\alpha_2^p=1\rangle, \\
\Phi_2(42)a2=& \langle \alpha,\alpha_1,\alpha_2:[\alpha_1,\alpha]=\alpha_1^p=\alpha_2,~\alpha^{p^4}=\alpha_2^p=1\rangle, \\
\Phi_2(411)b=& \langle \alpha,\alpha_1,\alpha_2,\gamma:[\alpha_1,\alpha]=\gamma^{p^3}=\alpha_2,~\alpha^p=\alpha_1^p=\alpha_2^p=1\rangle, \\
\Phi_2(411)c=& \langle \alpha,\alpha_1,\alpha_2:[\alpha_1,\alpha]=\alpha_2,~\alpha^{p^4}=\alpha_1^p=\alpha_2^p=1\rangle, \\
\Phi_2(33)=& \langle \alpha,\alpha_1,\alpha_2:[\alpha_1,\alpha]=\alpha_2=\alpha^{p^2},~\alpha_1^{p^3}=\alpha_2^p=1\rangle, \\
\Phi_2(321)c=& \langle \alpha,\alpha_1,\alpha_2,\gamma:[\alpha_1,\alpha]=\gamma^p=\alpha_2,~\alpha^{p^3}=\alpha_1^p=\alpha_2^p=1\rangle, \\
\Phi_2(321)d=& \langle \alpha,\alpha_1,\alpha_2,\gamma:[\alpha_1,\alpha]=\alpha_2=\gamma^{p^2},~\alpha^{p^2}=\alpha_1^p=\alpha_2^p=1\rangle, \\
\Phi_2(321)f=& \langle \alpha,\alpha_1,\alpha_2:[\alpha_1,\alpha]=\alpha_2,~\alpha^{p^3}=\alpha_1^{p^2}=\alpha_2^p=1\rangle, \\
\Phi_2(222)b=& \langle
\alpha,\alpha_1,\alpha_2,\gamma:[\alpha_1,\alpha]=\alpha_2=\gamma^p,~\alpha^{p^2}=\alpha_1^{p^2}=\alpha_2^p=1\rangle;
\end{align*}}
\item[{\rm (5)}]
{\allowdisplaybreaks
\begin{align*}
\Phi_5(3111)=& \langle \alpha_1,\alpha_2,\alpha_3,\alpha_4,\beta:[\alpha_1,\alpha_2]=[\alpha_3,\alpha_4]=\alpha_1^{p^2}=\beta,~\alpha_2^p=\alpha_3^p=\alpha_4^p=\beta^p=1\rangle, \\
\Phi_5(2211)a=& \langle \alpha_1,\alpha_2,\alpha_3,\alpha_4,\beta:[\alpha_1,\alpha_2]=[\alpha_3,\alpha_4]=\alpha_2^p=\beta,~\alpha_1^{p^2}=\alpha_3^p=\alpha_4^p=\beta^p=1\rangle, \\
\Phi_5(2211)b=& \langle \alpha_1,\alpha_2,\alpha_3,\alpha_4,\beta:[\alpha_1,\alpha_2]=[\alpha_3,\alpha_4]=\alpha_3^p=\beta,~\alpha_1^{p^2}=\alpha_2^p=\alpha_4^p=\beta^p=1\rangle, \\
\Phi_5(21^4)b=& \langle \alpha_1,\alpha_2,\alpha_3,\alpha_4,\beta,\gamma:[\alpha_1,\alpha_2]=[\alpha_3,\alpha_4]=\gamma^p=\beta,~\alpha_i^p=\beta^p=1~(i=1,2,3,4)\rangle, \\
\Phi_5(21^4)c=& \langle
\alpha_1,\alpha_2,\alpha_3,\alpha_4,\beta:[\alpha_1,\alpha_2]=[\alpha_3,\alpha_4]=\beta,~\alpha_1^{p^2}=\alpha_2^p=\alpha_3^p=\alpha_4^p=\beta^p=1\rangle.
\end{align*}}
\end{enumerate}

The pullbacks of order $p^6$ with an abelian quotient are precisely
those from families $(4),(12),(13)$ and $(15)$.

\begin{enumerate}
\item[{\rm (4)}]
{\allowdisplaybreaks
\begin{align*}
\Phi_4(321)a=& \langle \alpha,\alpha_1,\alpha_2,\beta_1,\beta_2:[\alpha_i,\alpha]=\beta_i,~\alpha^{p^2}=\beta_1,~\alpha_2^p=\beta_2,~ \alpha_1^p=\beta_i^p=1~ (i=1,2)\rangle, \\
\Phi_4(321)b=& \langle \alpha,\alpha_1,\alpha_2,\beta_1,\beta_2:[\alpha_i,\alpha]=\beta_i,~\alpha^{p^2}=\beta_1,~\alpha_1^p=\beta_2,~ \alpha_2^p=\beta_i^p=1~ (i=1,2)\rangle, \\
\Phi_4(321)c=& \langle \alpha,\alpha_1,\alpha_2,\beta_1,\beta_2:[\alpha_i,\alpha]=\beta_i,~\alpha^p=\beta_2,~\alpha_2^{p^2}=\beta_1,~ \alpha_1^p=\beta_i^p=1~ (i=1,2)\rangle, \\
\Phi_4(321)d=& \langle \alpha,\alpha_1,\alpha_2,\beta_1,\beta_2:[\alpha_i,\alpha]=\beta_i,~\alpha^p=\beta_2,~\alpha_1^{p^2}=\beta_1,~ \alpha_2^p=\beta_i^p=1~ (i=1,2)\rangle, \\
\Phi_4(321)e_r=& \langle \alpha,\alpha_1,\alpha_2,\beta_1,\beta_2:[\alpha_i,\alpha]=\beta_i,~\alpha_1^{p^2}=\beta_1,~ \alpha_2^p=\beta_2^r,~ \alpha^p=\beta_i^p=1~ (i=1,2)\rangle\\
 &\text{for}\ r=1,2,\dots,p-1, \\
\Phi_4(321)f_r=& \langle \alpha,\alpha_1,\alpha_2,\beta_1,\beta_2:[\alpha_i,\alpha]=\beta_i,~\alpha_1^p=\beta_2^r,~ \alpha_2^{p^2}=\beta_1,~ \alpha^p=\beta_i^p=1~ (i=1,2)\rangle\\
 &\text{for}\ r=1\ \text{or}\ \nu, \\
\Phi_4(3111)a=& \langle \alpha,\alpha_1,\alpha_2,\beta_1,\beta_2:[\alpha_i,\alpha]=\beta_i,~\alpha^{p^2}=\beta_1,~ \alpha_i^p=\beta_i^p=1~ (i=1,2)\rangle, \\
\Phi_4(3111)b=& \langle \alpha,\alpha_1,\alpha_2,\beta_1,\beta_2:[\alpha_i,\alpha]=\beta_i,~\alpha_1^{p^2}=\beta_1,~ \alpha^p=\alpha_2^p=\beta_i^p=1~ (i=1,2)\rangle, \\
\Phi_4(3111)c=& \langle \alpha,\alpha_1,\alpha_2,\beta_1,\beta_2:[\alpha_i,\alpha]=\beta_i,~\alpha_2^{p^2}=\beta_1,~ \alpha^p=\alpha_1^p=\beta_i^p=1~ (i=1,2)\rangle, \\
\Phi_4(222)a=& \langle \alpha,\alpha_1,\alpha_2,\beta_1,\beta_2:[\alpha_i,\alpha]=\beta_i=\alpha_i^p,~\alpha^{p^2}=\beta_i^p=1~ (i=1,2)\rangle, \\
\Phi_4(222)b_r=& \langle \alpha,\alpha_1,\alpha_2,\beta_1,\beta_2:[\alpha_i,\alpha]=\beta_i,~ \alpha_1^p=\beta_1^\kappa,~ \alpha_2^p=\beta_2,~\alpha^{p^2}=\beta_i^p=1~ (i=1,2)\rangle, \\
&\text{where}\ \kappa=g^r\ \text{for}\ r=1,2,\dots,\frac{1}{2}(p-1), \\
\Phi_4(222)c=& \langle \alpha,\alpha_1,\alpha_2,\beta_1,\beta_2:[\alpha_i,\alpha]=\beta_i,~ \alpha^p=\beta_1,~ \alpha_2^p=\beta_2,~\alpha_1^{p^2}=\beta_i^p=1~ (i=1,2)\rangle, \\
\Phi_4(222)d_1=& \langle \alpha,\alpha_1,\alpha_2,\beta_1,\beta_2:[\alpha_i,\alpha]=\beta_i,~ \alpha_1^p=\beta_2^{-1/4},~ \alpha_2^p=\beta_1\beta_2,~\alpha^{p^2}=\beta_i^p=1~ (i=1,2)\rangle, \\
\Phi_4(222)d_2=& \langle \alpha,\alpha_1,\alpha_2,\beta_1,\beta_2:[\alpha_i,\alpha]=\beta_i,~ \alpha^p=\beta_2,~ \alpha_2^p=\beta_1,~\alpha_1^{p^2}=\beta_i^p=1~ (i=1,2)\rangle, \\
\Phi_4(222)e_0=& \langle \alpha,\alpha_1,\alpha_2,\beta_1,\beta_2:[\alpha_i,\alpha]=\beta_i,~ \alpha_1^p=\beta_2,~ \alpha_2^p=\beta_1^\nu,~\alpha^{p^2}=\beta_i^p=1~ (i=1,2)\rangle, \\
\Phi_4(222)e_r=& \langle \alpha,\alpha_1,\alpha_2,\beta_1,\beta_2:[\alpha_i,\alpha]=\beta_i,~ \alpha_1^p=\beta_2^\kappa,~ \alpha_2^p=\beta_1\beta_2,~\alpha^{p^2}=\beta_i^p=1~ (i=1,2)\rangle, \\
&\text{where}\ 4\kappa=g^{2r+1}-1\ \text{for}\ r=1,2,\dots,\frac{1}{2}(p-1), \\
\Phi_4(2211)g=& \langle \alpha,\alpha_1,\alpha_2,\beta_1,\beta_2,\gamma:[\alpha_i,\alpha]=\beta_i,~ \gamma^p=\beta_2,~ \alpha^p=\beta_1,~ \alpha_i^p=\beta_i^p=1~ (i=1,2)\rangle, \\
\Phi_4(2211)h=& \langle \alpha,\alpha_1,\alpha_2,\beta_1,\beta_2,\gamma:[\alpha_i,\alpha]=\beta_i,~ \gamma^p=\beta_2,~ \alpha_1^p=\beta_1,~ \alpha^p=\alpha_2^p=\beta_i^p=1~ (i=1,2)\rangle, \\
\Phi_4(2211)i=& \langle \alpha,\alpha_1,\alpha_2,\beta_1,\beta_2,\gamma:[\alpha_i,\alpha]=\beta_i,~ \gamma^p=\beta_2,~ \alpha_2^p=\beta_1,~ \alpha^p=\alpha_1^p=\beta_i^p=1~ (i=1,2)\rangle, \\
\Phi_4(2211)j_1=& \langle \alpha,\alpha_1,\alpha_2,\beta_1,\beta_2:[\alpha_i,\alpha]=\beta_i,~ \alpha_1^p=\beta_1,~ \alpha^{p^2}=\alpha_2^p=\beta_i^p=1~ (i=1,2)\rangle, \\
\Phi_4(2211)j_2=& \langle \alpha,\alpha_1,\alpha_2,\beta_1,\beta_2:[\alpha_i,\alpha]=\beta_i,~ \alpha^p=\beta_1,~ \alpha_1^{p^2}=\alpha_2^p=\beta_i^p=1~ (i=1,2)\rangle, \\
\Phi_4(2211)k=& \langle \alpha,\alpha_1,\alpha_2,\beta_1,\beta_2:[\alpha_i,\alpha]=\beta_i,~ \alpha^p=\beta_2,~ \alpha_1^{p^2}=\alpha_2^p=\beta_i^p=1~ (i=1,2)\rangle, \\
\Phi_4(2211)l=& \langle \alpha,\alpha_1,\alpha_2,\beta_1,\beta_2:[\alpha_i,\alpha]=\beta_i,~ \alpha_2^p=\beta_1,~ \alpha^{p^2}=\alpha_1^p=\beta_i^p=1~ (i=1,2)\rangle, \\
\Phi_4(2211)m=& \langle \alpha,\alpha_1,\alpha_2,\beta_1,\beta_2:[\alpha_i,\alpha]=\beta_i,~ \alpha_2^p=\beta_2,~ \alpha_1^{p^2}=\alpha^p=\beta_i^p=1~ (i=1,2)\rangle, \\
\Phi_4(2211)n=& \langle \alpha,\alpha_1,\alpha_2,\beta_1,\beta_2:[\alpha_i,\alpha]=\beta_i,~ \alpha_2^p=\beta_1,~ \alpha_1^{p^2}=\alpha^p=\beta_i^p=1~ (i=1,2)\rangle, \\
\Phi_4(21^4)d=& \langle \alpha,\alpha_1,\alpha_2,\beta_1,\beta_2,\gamma:[\alpha_i,\alpha]=\beta_i,~ \gamma^p=\beta_1,~ \alpha^p=\alpha_i^p=\beta_i^p=1~ (i=1,2)\rangle, \\
\Phi_4(21^4)e=& \langle \alpha,\alpha_1,\alpha_2,\beta_1,\beta_2:[\alpha_i,\alpha]=\beta_i,~ \alpha^{p^2}=\alpha_i^p=\beta_i^p=1~ (i=1,2)\rangle, \\
\Phi_4(21^4)f=& \langle
\alpha,\alpha_1,\alpha_2,\beta_1,\beta_2:[\alpha_i,\alpha]=\beta_i,~
\alpha_1^{p^2}=\alpha_2^p=\alpha^p=\beta_i^p=1~ (i=1,2)\rangle;
\end{align*}}
\item[{\rm (12)}]
{\allowdisplaybreaks
\begin{align*}
\Phi_{12}(2211)a=& \langle
\alpha_1,\alpha_2,\beta_1,\beta_2,\gamma_1,\gamma_2:[\alpha_i,\beta_i]=\gamma_i,~
\alpha_1^p=\gamma_1,~ \beta_1^p=\gamma_2,~ \alpha_2^p=\beta_2^p=\gamma_i^p=1~ (i=1,2)\rangle,\\
\Phi_{12}(2211)c=& \langle
\alpha_1,\alpha_2,\beta_1,\beta_2,\gamma_1,\gamma_2:[\alpha_i,\beta_i]=\gamma_i,~
\alpha_1^p=\gamma_1\gamma_2,~ \alpha_2^p=\gamma_2,~ \beta_i^p=\gamma_i^p=1~ (i=1,2)\rangle,\\
\Phi_{12}(2211)d=& \langle
\alpha_1,\alpha_2,\beta_1,\beta_2,\gamma_1,\gamma_2:[\alpha_i,\beta_i]=\gamma_i,~
\alpha_1^p=\gamma_2,~ \alpha_2^p=\gamma_1,~ \beta_i^p=\gamma_i^p=1~ (i=1,2)\rangle,\\
\Phi_{12}(2211)e=& \langle
\alpha_1,\alpha_2,\beta_1,\beta_2,\gamma_1,\gamma_2:[\alpha_i,\beta_i]=\gamma_i,~
\alpha_1^p=\gamma_1\gamma_2,~ \alpha_2^p=\gamma_1,~ \beta_i^p=\gamma_i^p=1~ (i=1,2)\rangle,\\
\Phi_{12}(2211)f=& \langle
\alpha_1,\alpha_2,\beta_1,\beta_2,\gamma_1,\gamma_2:[\alpha_i,\beta_i]=\gamma_i,~
\alpha_1^p=\alpha_2^p=\gamma_1,~ \beta_1^p=\gamma_2,~ \beta_2^p=\gamma_i^p=1~ (i=1,2)\rangle,\\
\Phi_{12}(2211)g=& \langle
\alpha_1,\alpha_2,\beta_1,\beta_2,\gamma_1,\gamma_2:[\alpha_i,\beta_i]=\gamma_i,~
\alpha_1^p=\gamma_1,~ \alpha_2^p=\beta_1^p=\gamma_2,~ \beta_2^p=\gamma_i^p=1~ (i=1,2)\rangle,\\
\Phi_{12}(2211)h=& \langle
\alpha_1,\alpha_2,\beta_1,\beta_2,\gamma_1,\gamma_2:[\alpha_i,\beta_i]=\gamma_i,~
\alpha_1^p=\beta_2^p=\gamma_1,~ \alpha_2^p=\beta_1^p=\gamma_2,~ \gamma_i^p=1~ (i=1,2)\rangle,\\
\Phi_{12}(2211)i=& \langle
\alpha_1,\alpha_2,\beta_1,\beta_2,\gamma_1,\gamma_2:[\alpha_i,\beta_i]=\gamma_i,~
\alpha_1^p=\gamma_1,~ \alpha_2^p=\gamma_1\gamma_2,~ \beta_1^p=\gamma_2,~ \beta_2^p=\gamma_i^p=1\\
& (i=1,2)\rangle,\\
\Phi_{12}(21^4)b=& \langle
\alpha_1,\alpha_2,\beta_1,\beta_2,\gamma_1,\gamma_2:[\alpha_i,\beta_i]=\gamma_i,~
\alpha_1^p=\gamma_1\gamma_2,~ \alpha_2^p=\beta_i^p=\gamma_i^p=1~ (i=1,2)\rangle,\\
\Phi_{12}(21^4)c=& \langle
\alpha_1,\alpha_2,\beta_1,\beta_2,\gamma_1,\gamma_2:[\alpha_i,\beta_i]=\gamma_i,~
\alpha_1^p=\gamma_2,~ \alpha_2^p=\beta_i^p=\gamma_i^p=1~ (i=1,2)\rangle,\\
\Phi_{12}(21^4)d=& \langle
\alpha_1,\alpha_2,\beta_1,\beta_2,\gamma_1,\gamma_2:[\alpha_i,\beta_i]=\gamma_i,~
\alpha_1^p=\alpha_2^p=\gamma_1,~ \beta_i^p=\gamma_i^p=1~ (i=1,2)\rangle,\\
\Phi_{12}(21^4)e=& \langle
\alpha_1,\alpha_2,\beta_1,\beta_2,\gamma_1,\gamma_2:[\alpha_i,\beta_i]=\gamma_i,~
\alpha_1^p=\alpha_2^p=\gamma_1\gamma_2,~ \beta_i^p=\gamma_i^p=1~
(i=1,2)\rangle;
\end{align*}}
\item[{\rm (13)}]
{\allowdisplaybreaks
\begin{align*}
\Phi_{13}(2211)a=& \langle
\alpha_1,\dots,\alpha_4,\beta_1,\beta_2:[\alpha_1,\alpha_{i+1}]=\beta_i,~ [\alpha_2,\alpha_4]=\alpha_2^p=
\beta_2,~ \alpha_1^p=\beta_1,\\
&\alpha_3^p=\alpha_4^p=\beta_i^p=1~ (i=1,2)\rangle,\\
\Phi_{13}(2211)b=& \langle
\alpha_1,\dots,\alpha_4,\beta_1,\beta_2:[\alpha_1,\alpha_{i+1}]=\beta_i,~
[\alpha_2,\alpha_4]=\alpha_3^p=
\beta_2,~ \alpha_1^p=\beta_1,\\
&\alpha_2^p=\alpha_4^p=\beta_i^p=1~ (i=1,2)\rangle,\\
\Phi_{13}(2211)c_r=& \langle
\alpha_1,\dots,\alpha_4,\beta_1,\beta_2:[\alpha_1,\alpha_{i+1}]=\beta_i,~
[\alpha_2,\alpha_4]^r=\alpha_2^p=
\beta_2^r,~ \alpha_3^p=\beta_1,\\
&\alpha_1^p=\alpha_4^p=\beta_i^p=1~ (i=1,2)\rangle\ \text{for}\ r=1\ \text{or}\ \nu,\\
\Phi_{13}(2211)d=& \langle
\alpha_1,\dots,\alpha_4,\beta_1,\beta_2:[\alpha_1,\alpha_{i+1}]=\beta_i,~
[\alpha_2,\alpha_4]=\alpha_1^p=
\beta_2,~ \alpha_3^p=\beta_1,\\
&\alpha_2^p=\alpha_4^p=\beta_i^p=1~ (i=1,2)\rangle,\\
\Phi_{13}(2211)e_r=& \langle
\alpha_1,\dots,\alpha_4,\beta_1,\beta_2:[\alpha_1,\alpha_{i+1}]=\beta_i,~
[\alpha_2,\alpha_4]^r=\alpha_4^p=
\beta_2^r,~ \alpha_1^p=\beta_1,\\
&\alpha_2^p=\alpha_3^p=\beta_i^p=1~ (i=1,2)\rangle\ \text{for}\ r=1,\dots, p-1,\\
\Phi_{13}(2211)f=& \langle
\alpha_1,\dots,\alpha_4,\beta_1,\beta_2:[\alpha_1,\alpha_{i+1}]=\beta_i,~
[\alpha_2,\alpha_4]=\alpha_4^p=
\beta_2,~ \alpha_3^p=\beta_1,\\
&\alpha_1^p=\alpha_2^p=\beta_i^p=1~ (i=1,2)\rangle,\\
\Phi_{13}(21^4)a=& \langle
\alpha_1,\dots,\alpha_4,\beta_1,\beta_2:[\alpha_1,\alpha_{i+1}]=\beta_i,~
[\alpha_2,\alpha_4]=\beta_2,~ \alpha_1^p=\beta_1,\\
&\alpha_{i+1}^p=\alpha_4^p=\beta_i^p=1~ (i=1,2)\rangle,\\
\Phi_{13}(21^4)b=& \langle
\alpha_1,\dots,\alpha_4,\beta_1,\beta_2:[\alpha_1,\alpha_{i+1}]=\beta_i,~
[\alpha_2,\alpha_4]=\alpha_1^p=\beta_2,\\
&\alpha_{i+1}^p=\alpha_4^p=\beta_i^p=1~ (i=1,2)\rangle,\\
\Phi_{13}(21^4)c=& \langle
\alpha_1,\dots,\alpha_4,\beta_1,\beta_2:[\alpha_1,\alpha_{i+1}]=\beta_i,~
[\alpha_2,\alpha_4]=\alpha_3^p=\beta_2,\\
&\alpha_i^p=\alpha_4^p=\beta_i^p=1~ (i=1,2)\rangle,\\
\Phi_{13}(21^4)d=& \langle
\alpha_1,\dots,\alpha_4,\beta_1,\beta_2:[\alpha_1,\alpha_{i+1}]=\beta_i,~
[\alpha_2,\alpha_4]=\beta_2,~ \alpha_3^p=\beta_1,\\
&\alpha_i^p=\alpha_4^p=\beta_i^p=1~ (i=1,2)\rangle,\\
\Phi_{13}(1^6)=& \langle
\alpha_1,\dots,\alpha_4,\beta_1,\beta_2:[\alpha_1,\alpha_{i+1}]=\beta_i,~
[\alpha_2,\alpha_4]=\beta_2,\\
&\alpha_i^p=\alpha_3^p=\alpha_4^p=\beta_i^p=1~ (i=1,2)\rangle;
\end{align*}}
\item[{\rm (15)}]
{\allowdisplaybreaks
\begin{align*}
\Phi_{15}(2211)a=& \langle
\alpha_1,\dots,\alpha_4,\beta_1,\beta_2:[\alpha_1,\alpha_{i+1}]=\beta_i,~
[\alpha_3,\alpha_4]=\alpha_1^p=\beta_1,~ [\alpha_2,\alpha_4]=\alpha_2^{gp}=\beta_2^g,\\
&\alpha_3^p=\alpha_4^p=\beta_i^p=1~ (i=1,2)\rangle,\\
\Phi_{15}(2211)b_{r,s}=& \langle
\alpha_1,\dots,\alpha_4,\beta_1,\beta_2:[\alpha_1,\alpha_{i+1}]=\beta_i,~
[\alpha_3,\alpha_4]=\beta_1,~ [\alpha_2,\alpha_4]^\kappa=\alpha_2^{gp}=\beta_2^{g\kappa},\\
&\alpha_1^p=\beta_1\beta_2^r,~ \alpha_3^p=\alpha_4^p=\beta_i^p=1~
(i=1,2)\rangle,\ \text{where}\ \kappa=g^{[1/(2n)]+s}\ \text{and}\\
&g^n=g^2(g-r^2)\ \text{for}\ r=1,2,\dots,\frac{1}{2}(p-1)\ \text{and}\ s=0,1,\dots,m,\ \text{with}\\
&m=\frac{1}{2}(p-3)+n-2[1/(2n)]\ \text{and}\ [1/(2n)]=\text{integral
part of}\ 1/(2n),\\
\Phi_{15}(2211)c=& \langle
\alpha_1,\dots,\alpha_4,\beta_1,\beta_2:[\alpha_1,\alpha_{i+1}]=\beta_i,~
[\alpha_3,\alpha_4]=\alpha_1^p=\beta_1,~ [\alpha_2,\alpha_4]=\alpha_4^{-p}=\beta_2^g,\\
&\alpha_{i+1}^p=\beta_i^p=1~ (i=1,2)\rangle,\\
\Phi_{15}(2211)d_r=& \langle
\alpha_1,\dots,\alpha_4,\beta_1,\beta_2:[\alpha_1,\alpha_{i+1}]=\beta_i,~
[\alpha_3,\alpha_4]=\alpha_1^p=\beta_1,~ [\alpha_2,\alpha_4]=\beta_2^g,\\
&\alpha_4^p=\beta_2^\kappa,~ \alpha_{i+1}^p=\beta_i^p=1~ (i=1,2)\rangle,\ \text{where}\ \kappa=g^r\ \text{for}\ r=1,2,\dots,\frac{1}{2}(p-1),\\
\Phi_{15}(21^4)=& \langle
\alpha_1,\dots,\alpha_4,\beta_1,\beta_2:[\alpha_1,\alpha_{i+1}]=\beta_i,~
[\alpha_3,\alpha_4]=\alpha_1^p=\beta_1,~ [\alpha_2,\alpha_4]=\beta_2^g,\\
&\alpha_{i+1}^p=\alpha_4^p=\beta_i^p=1~ (i=1,2)\rangle,\\
\Phi_{15}(1^6)=& \langle
\alpha_1,\dots,\alpha_4,\beta_1,\beta_2:[\alpha_1,\alpha_{i+1}]=\beta_i,~
[\alpha_3,\alpha_4]=\beta_1,~ [\alpha_2,\alpha_4]=\beta_2^g,\\
&\alpha_i^p=\alpha_3^p=\alpha_4^p=\beta_i^p=1~ (i=1,2)\rangle.
\end{align*}}
\end{enumerate}

Finally, the groups of order $p^6$ having quotient $(C_{p^2})^2$
(obtained by factoring out a central cyclic group of order $p^2$)
are precisely the groups from family $(14)$.

\begin{enumerate}
\item[{\rm (14)}]
{\allowdisplaybreaks
\begin{align*}
\Phi_{14}(42)=& \langle \alpha_1,\alpha_2,\beta:[\alpha_1,\alpha_2]=\beta,~\alpha_1^{p^2}=\beta,~\alpha_2^{p^2}=\beta^{p^2}=1\rangle, \\
\Phi_{14}(321)=& \langle \alpha_1,\alpha_2,\beta:[\alpha_1,\alpha_2]=\beta,~\alpha_1^{p^2}=\beta^p,~\alpha_2^{p^2}=\beta^{p^2}=1\rangle, \\
\Phi_{14}(222)=& \langle
\alpha_1,\alpha_2,\beta:[\alpha_1,\alpha_2]=\beta,~\alpha_1^{p^2}=\alpha_2^{p^2}=\beta^{p^2}=1\rangle.
\end{align*}}
\end{enumerate}

%--------------------------------------Section 4-------------------------------------------------
\section{Results and sample proofs}
\label{4}

Our results are displayed in Tables $1,\dots,6$, which give
necessary and sufficient conditions for the realizability of a group
over a field $k$. The tables are organized by the criteria that we
use for solving the embedding problems -- Theorem \ref{c3.6},
Theorem \ref{t4.0} and Theorem \ref{t3.7}. In each table, the first
column indicates the group, according to the presentations given in
Section \ref{3}. The second column gives labels for the elements of
$k^\times$ that are required to be independent$\pmod{k^{\times p}}$.
The third column indicates the assumption that certain root of unity
is in $k$. Note that this is not a necessary condition in general,
but we need these roots of unity so that we can apply Theorem
\ref{c3.6}. The fourth column gives the obstructions to the
realizability of these groups, i.e., the $p$-cyclic algebras
required to be trivial in the Brauer group $\Br(k)$.

We give detailed proofs for two representative groups. Proofs for
the other groups follow the same general outline. The first group is
representative for the groups of orders $p^5$ and $p^6$ from
families $(2)$ and $(5)$.

\newtheorem{t4.1}{Theorem}[section]
\begin{t4.1}\label{t4.1}
Let $k$ be a field of characteristic not $p$ containing a primitive
$p^3$-th root of unity $\zeta_{p^3}$. The group $\Phi_2(41)$ is
realizable as a Galois group over $k$ if and only if there exist
elements $a_1,a_2\in k^\times$ independent$\pmod{k^{\times p}}$ such
that $(\zeta_{p^3}^{-1}a_1,a_2;\zeta)=1\in\Br(k)$.
\end{t4.1}
\begin{proof}
Since $\zeta_{p^3}\in k$, the group $C_p\times C_{p^3}$ is
realizable over $k$ if and only if there exist elements $a_1,a_2\in
k^\times$ independent$\pmod{k^{\times p}}$. Then
$K/k=k(\root{p}\of{a_1},\root{p^3}\of{a_2})/k$ is a $C_p\times
C_{p^3}$ extension, and we may consider the embedding problem
$(K/k,\Phi_2(41),\mu_p)$ related to the group extension
\begin{equation*}
1\longrightarrow \langle\alpha_2\rangle\cong \mu_p\longrightarrow
\Phi_2(41)\underset{\attop{\alpha\mapsto\sigma_2}{
\alpha_1\mapsto\sigma_1}}{\longrightarrow} C_p\times
C_{p^3}\longrightarrow 1,
\end{equation*}
where $C_p\times
C_{p^3}=\langle\sigma_1\rangle\times\langle\sigma_2\rangle$,
according to the notation of Theorem \ref{c3.6}. Let $s_1=\alpha_1$
and $s_2=\alpha$ in $\Phi_2(41)$. Note that
$t=2,n_1=1,n_2=3,m_1=0,m_2=1,r=2,n=3$ and $d_{21}=1$. Therefore, the
obstruction to the embedding problem $(K/k,\Phi_2(41),\mu_p)$ is
$(a_1,a_2;\zeta)(a_2,\zeta_{p^3};\zeta)=(\zeta_{p^3}^{-1}a_1,a_2;\zeta)\in\Br(k)$.
\end{proof}

The second group is representative for all pullbacks of orders $p^5$
and $p^6$ that have an abelian quotient.

\newtheorem{t4.2}[t4.1]{Theorem}
\begin{t4.2}\label{t4.2}
Let $k$ be a field of characteristic not $p$ containing a primitive
$p$-th root of unity $\zeta$. The group $\Phi_4(221)a$ is realizable
as a Galois group over $k$ if and only if there exist elements
$a_1,a_2,a_3\in k^\times$ independent$\pmod{k^{\times p}}$ such that
$(\zeta^{-1} a_2,a_3;\zeta)=(a_1,\zeta a_3;\zeta)=1\in\Br(k)$.
\end{t4.2}
\begin{proof}
Since $\zeta\in k$, the group $(C_p)^3$ is realizable over $k$ if
and only if there exist elements $a_1,a_2,a_3\in k^\times$ that are
independent$\pmod{k^{\times p}}$. Assume that such elements exist
and let
$K/k=k(\root{p}\of{a_1},\root{p}\of{a_2},\root{p}\of{a_3})/k$. Let
$\sigma_i, i=1,2,3$ be the generators of  $(C_p)^3$ which act on
$K/k$ according to the description given in Theorem \ref{c3.6}.

Next, using the notation of Theorem \ref{t4.0} put
$N_i=\langle\beta_i\rangle$ for $i=1,2$. We have the group extension
\begin{equation}\label{e0}
1\longrightarrow N_1N_2\cong (\mu_p)^2\longrightarrow
\Phi_4(221)a\underset{\scriptsize\begin{array}{c}
               \alpha_1 \mapsto \sigma_1\\
               \alpha_2 \mapsto \sigma_2\\
               \alpha\mapsto \sigma_3\\
                \end{array}}{\longrightarrow}  (C_p)^3\longrightarrow
                1.
\end{equation}
By Theorem \ref{t4.0}, the embedding problem given by $K/k$ and
\eqref{e0} is solvable if and only if the embedding problems given
by $K/k$ and the group extensions
\begin{equation}\label{e1}
1\longrightarrow N_2\cong \mu_p\longrightarrow
\Phi_4(221)a/N_1\underset{\scriptsize\begin{array}{c}
               \alpha_1 \mapsto \sigma_1\\
               \alpha_2 \mapsto \sigma_2\\
               \alpha\mapsto \sigma_3\\
                \end{array}}{\longrightarrow}  (C_p)^3\longrightarrow
                1
\end{equation}
and
\begin{equation}\label{e2}
1\longrightarrow N_1\cong \mu_p\longrightarrow
\Phi_4(221)a/N_2\underset{\scriptsize\begin{array}{c}
               \alpha_1 \mapsto \sigma_1\\
               \alpha_2 \mapsto \sigma_2\\
               \alpha\mapsto \sigma_3\\
                \end{array}}{\longrightarrow}  (C_p)^3\longrightarrow
                1
\end{equation}
are solvable. In the notation of Theorem \ref{c3.6}, let
$s_1=\alpha_1,s_2=\alpha_2$ and $s_3=\alpha$. Note that for
\eqref{e1} we have $t=3,r=3,m_3=1,d_{32}=1$, so the obstruction is
$(a_2,a_3;\zeta)(a_3,\zeta;\zeta)=(\zeta^{-1} a_2,a_3;\zeta)$.  For
\eqref{e2} we have $t=3,r=1,m_1=1, d_{31}=1$, so the obstruction is
$(a_1,a_3;\zeta)(a_1,\zeta;\zeta)=(a_1,\zeta a_3;\zeta)$.  We are
done.
\end{proof}

For the three groups of order $p^6$ that have quotient $(C_{p^2})^2$
obtained by factoring out $\mu_{p^2}$, we have only to modify the
proof of Theorem \ref{t4.1}, so that we apply Theorem \ref{t3.7}
instead of Theorem \ref{c3.6}. Note that the kernel $\langle
[\alpha_1,\alpha_2]\rangle$ is contained in the Frattini subgroup
$\Phi(G)=[G,G]\cdot G^p$, so the obstructions are again in terms of
proper solvability.

{\bf Standing notations.} For the groups of orders $p^5$ and $p^6$
from family $(2)$ we put $s_3=\gamma$ (where appears the generator
$\gamma$). For the groups of order $p^6$ from family $(12)$ we put
$s_1=\alpha_1,s_2=\alpha_2,s_3=\beta_1,s_4=\beta_2$. In general, if
$C_{p^s}$ is a direct factor of a given abelian quotient we suppose
that $k$ contains a primitive $p^s$-th root of unity so that
$C_{p^s}$ is always realizable. On some occasions, however, we can
lessen the condition for the roots of unity and suppose only that a
primitive $p^{s-1}$-th root of unity is contained in $k$. In this
case, we have to require that $C_{p^s}$ is realizable over $k$. That
is why in the tables sometimes appears additional obstructions of
the kind $(*,\zeta_{p^{s-1}};\zeta)$.

\bigskip

\begin{table}[ht]
\caption{Groups of order $p^5$ with abelian quotients}\centering
\begin{tabular}{|l|l|l|l|} \hline
Group & $\vert k^\times/k^{\times p}\vert$ & Root of unity & Trivial elements in \Br(k)   \\
\hline\hline $\Phi_2(41)$ & $a_1,a_2$ & $\zeta_{p^3}\in k$ &
$(\zeta_{p^3}^{-1}a_1,a_2;\zeta)$\\
\hline $\Phi_2(32)a1$ & $a_1,a_2$ & $\zeta_{p^2}\in k$ &
$(\zeta_{p^2}^{-1}a_1,a_2;\zeta)$\\
\hline $\Phi_2(32)a2$ & $a_1,a_2$ & $\zeta_{p^2}\in k$ &
$(a_2,\zeta_{p^2};\zeta),~ (a_1,a_2;\zeta)$\\
\hline $\Phi_2(311)b$ & $a_1,a_2,a_3$ & $\zeta_{p^2}\in k$ &
$(a_1,a_2;\zeta)(a_3,\zeta_{p^2};\zeta)$\\
\hline $\Phi_2(311)c$ & $a_1,a_2$ & $\zeta_{p^2}\in k$ &
$(a_2,\zeta_{p^2};\zeta),~ (a_1,a_2;\zeta)$\\
\hline $\Phi_2(221)c$ & $a_1,a_2,a_3$ & $\zeta\in k$ &
$(a_2,\zeta;\zeta),~ (a_1,a_2;\zeta)(a_3,\zeta;\zeta)$\\
\hline $\Phi_2(221)d$ & $a_1,a_2$ & $\zeta\in k$ &
$(a_1,\zeta;\zeta),~ (a_2,\zeta;\zeta),~ (a_1,a_2;\zeta)$\\
\hline $\Phi_5(2111)$ & $a_1,a_2,a_3,a_4$ & $\zeta\in k$ &
$(a_1,\zeta a_2;\zeta)(a_3,a_4;\zeta)$\\
\hline $\Phi_5(1^5)$ & $a_1,a_2,a_3,a_4$ & $\zeta\in k$ &
$(a_1,a_2;\zeta)(a_3,a_4;\zeta)$\\
\hline
\end{tabular}
\end{table}

\clearpage

\begin{table}[ht]
\caption{Pullbacks of order $p^5$ with abelian quotients}\centering
\begin{tabular}{|l|l|l|l|} \hline
Group & $\vert k^\times/k^{\times p}\vert$ & Root of unity & Trivial elements in \Br(k)   \\
\hline\hline $\Phi_4(221)a$ & $a_1,a_2,a_3$ & $\zeta\in k$ &
$(\zeta^{-1} a_2,a_3;\zeta),~ (a_1,\zeta a_3;\zeta)$\\
\hline $\Phi_4(221)b$ & $a_1,a_2,a_3$ & $\zeta\in k$ &
$(\zeta^{-1} a_2,a_3;\zeta),~ (a_1,a_3;\zeta)(a_2,\zeta;\zeta)$\\
\hline $\Phi_4(221)c$ & $a_1,a_2,a_3$ & $\zeta\in k$ &
$(a_2,\zeta a_3;\zeta),~ (a_1,\zeta a_3;\zeta)$\\
\hline $\Phi_4(221)d_r$ & $a_1,a_2,a_3$ & $\zeta\in k$ &
$(a_2,\zeta a_3;\zeta),~ (a_1,\zeta^\kappa a_3;\zeta)$\\
\hline $\Phi_4(221)e$ & $a_1,a_2,a_3$ & $\zeta\in k$ &
$(a_2,\zeta a_3;\zeta)(a_1,\zeta^{-1/4};\zeta),~ (a_1,a_3;\zeta)(a_2,\zeta;\zeta)$\\
\hline $\Phi_4(221)f_0$ & $a_1,a_2,a_3$ & $\zeta\in k$ &
$(a_2,a_3;\zeta)(a_1,\zeta;\zeta),~ (a_1,a_3;\zeta)(a_2,\zeta^\nu;\zeta)$\\
\hline $\Phi_4(221)f_r$ & $a_1,a_2,a_3$ & $\zeta\in k$ &
$(a_2,\zeta a_3;\zeta)(a_1,\zeta^\kappa;\zeta),~ (a_1,a_3;\zeta)(a_2,\zeta;\zeta)$\\
\hline $\Phi_4(2111)a$ & $a_1,a_2,a_3$ & $\zeta\in k$ &
$(\zeta^{-1} a_2,a_3;\zeta),~ (a_1,a_3;\zeta)$\\
\hline $\Phi_4(2111)b$ & $a_1,a_2,a_3$ & $\zeta\in k$ &
$(a_2,a_3;\zeta),~ (a_1,\zeta a_3;\zeta)$\\
\hline $\Phi_4(2111)c$ & $a_1,a_2,a_3$ & $\zeta\in k$ &
$(a_2,a_3;\zeta),~ (a_1,a_3;\zeta)(a_2,\zeta;\zeta)$\\
\hline $\Phi_4(1^5)$ & $a_1,a_2,a_3$ & $\zeta\in k$ &
$(a_2,a_3;\zeta),~ (a_1,a_3;\zeta)$\\
\hline
\end{tabular}
\end{table}

\bigskip

\begin{table}[ht]
\caption{Groups of order $p^6$ with abelian quotients}\centering
\begin{tabular}{|l|l|l|l|} \hline
Group & $\vert k^\times/k^{\times p}\vert$ & Root of unity & Trivial elements in \Br(k)   \\
\hline\hline $\Phi_2(51)$ & $a_1,a_2$ & $\zeta_{p^4}\in k$ &
$(\zeta_{p^4}^{-1}a_1,a_2;\zeta)$\\
\hline $\Phi_2(42)a1$ & $a_1,a_2$ & $\zeta_{p^3}\in k$ &
$(\zeta_{p^3}^{-1}a_1,a_2;\zeta)$\\
\hline $\Phi_2(42)a2$ & $a_1,a_2$ & $\zeta_{p^3}\in k$ &
$(a_2,\zeta_{p^3};\zeta),~ (a_1,a_2;\zeta)$\\
\hline $\Phi_2(411)b$ & $a_1,a_2,a_3$ & $\zeta_{p^3}\in k$ &
$(a_1,a_2;\zeta)(a_3,\zeta_{p^3};\zeta)$\\
\hline $\Phi_2(411)c$ & $a_1,a_2$ & $\zeta_{p^3}\in k$ &
$(a_2,\zeta_{p^3};\zeta),~ (a_1,a_2;\zeta)$\\
\hline $\Phi_2(33)$ & $a_1,a_2$ & $\zeta_{p^2}\in k$ &
$(a_1,\zeta_{p^2};\zeta),~ (\zeta_{p^2}^{-1} a_1,a_2;\zeta)$\\
\hline $\Phi_2(321)c$ & $a_1,a_2,a_3$ & $\zeta_{p^2}\in k$ &
$(a_2,\zeta_{p^2};\zeta),~ (a_1,a_2;\zeta)$\\
\hline $\Phi_2(321)d$ & $a_1,a_2,a_3$ & $\zeta_{p^2}\in k$ &
$(a_1,a_2;\zeta)(a_3,\zeta_{p^2};\zeta)$\\
\hline $\Phi_2(321)f$ & $a_1,a_2$ & $\zeta_{p^2}\in k$ &
$(a_2,\zeta_{p^2};\zeta),~(a_1,a_2;\zeta)$\\
\hline $\Phi_2(222)b$ & $a_1,a_2,a_3$ & $\zeta\in k$ &
$(a_1,\zeta;\zeta),~(a_2,\zeta;\zeta),~ (a_1,a_2;\zeta)(a_3,\zeta;\zeta)$\\
\hline $\Phi_5(3111)$ & $a_1,a_2,a_3,a_4$ & $\zeta_{p^2}\in k$ &
$(a_1,\zeta_{p^2} a_2;\zeta)(a_3,a_4;\zeta)$\\
\hline $\Phi_5(2211)a$ & $a_1,a_2,a_3,a_4$ & $\zeta\in k$ &
$(a_1,\zeta;\zeta),~ (\zeta^{-1} a_1, a_2;\zeta)(a_3,a_4;\zeta)$\\
\hline $\Phi_5(2211)b$ & $a_1,a_2,a_3,a_4$ & $\zeta\in k$ &
$(a_1,\zeta;\zeta),~ (a_1, a_2;\zeta)(a_3,\zeta a_4;\zeta)$\\
\hline $\Phi_5(21^4)b$ & $a_1,a_2,a_3,a_4,a_5$ & $\zeta\in k$ &
$(a_1, a_2;\zeta)(a_3,a_4;\zeta)(a_5,\zeta;\zeta)$\\
\hline $\Phi_5(21^4)c$ & $a_1,a_2,a_3,a_4$ & $\zeta\in k$ &
$(a_1,\zeta;\zeta),~ (a_1, a_2;\zeta)(a_3,a_4;\zeta)$\\
\hline
\end{tabular}
\end{table}

\clearpage

\begin{table}[ht]
\caption{Pullbacks of order $p^6$ with abelian quotients: Part
I}\centering
\begin{tabular}{|l|l|l|l|} \hline
Group & $\vert k^\times/k^{\times p}\vert$ & Root of unity & Trivial elements in \Br(k)   \\
\hline\hline $\Phi_4(321)a$ & $a_1,a_2,a_3$ & $\zeta_{p^2}\in k$ &
$(a_2,a_3;\zeta),~ (\zeta_{p^2}^{-1} a_1,a_3;\zeta)$\\
\hline $\Phi_4(321)b$ & $a_1,a_2,a_3$ & $\zeta_{p^2}\in k$ &
$(a_2,a_3;\zeta),~ (\zeta_{p^2}^{-1} a_1,a_3;\zeta)$\\
\hline $\Phi_4(321)c$ & $a_1,a_2,a_3$ & $\zeta_{p^2}\in k$ &
$(a_2,a_3;\zeta),~ (a_1,a_3;\zeta)(a_2,\zeta_{p^2};\zeta)$\\
\hline $\Phi_4(321)d$ & $a_1,a_2,a_3$ & $\zeta_{p^2}\in k$ &
$(a_2,a_3;\zeta),~ (a_1,\zeta_{p^2} a_3;\zeta)$\\
\hline $\Phi_4(321)e_r$ & $a_1,a_2,a_3$ & $\zeta_{p^2}\in k$ &
$(a_2,a_3;\zeta),~ (a_1,\zeta_{p^2} a_3;\zeta)$\\
\hline $\Phi_4(321)f_r$ & $a_1,a_2,a_3$ & $\zeta_{p^2}\in k$ &
$(a_2,a_3;\zeta),~ (a_1,a_3;\zeta)(a_2,\zeta_{p^2};\zeta)$\\
\hline $\Phi_4(3111)a$ & $a_1,a_2,a_3$ & $\zeta_{p^2}\in k$ &
$(a_2,a_3;\zeta),~ (\zeta_{p^2}^{-1} a_1,a_3;\zeta)$\\
\hline $\Phi_4(3111)b$ & $a_1,a_2,a_3$ & $\zeta_{p^2}\in k$ &
$(a_2,a_3;\zeta),~ (a_1,\zeta_{p^2} a_3;\zeta)$\\
\hline $\Phi_4(3111)c$ & $a_1,a_2,a_3$ & $\zeta_{p^2}\in k$ &
$(a_2,a_3;\zeta),~ (a_1,a_3;\zeta)(a_2,\zeta_{p^2};\zeta)$\\
\hline $\Phi_4(222)a$ & $a_1,a_2,a_3$ & $\zeta\in k$ &
$(a_3,\zeta;\zeta),~ (a_2,\zeta a_3;\zeta),~ (a_1,\zeta a_3;\zeta)$\\
\hline $\Phi_4(222)b_r$ & $a_1,a_2,a_3$ & $\zeta\in k$ &
$(a_3,\zeta;\zeta),~ (a_2,\zeta a_3;\zeta),~ (a_1,\zeta^\kappa a_3;\zeta)$\\
\hline $\Phi_4(222)c$ & $a_1,a_2,a_3$ & $\zeta\in k$ &
$(a_1,\zeta;\zeta),~ (a_2,\zeta a_3;\zeta),~ (\zeta^{-1} a_1, a_3;\zeta)$\\
\hline $\Phi_4(222)d_1$ & $a_1,a_2,a_3$ & $\zeta\in k$ &
$(a_3,\zeta;\zeta),~ (a_2,\zeta a_3;\zeta)(a_1,\zeta^{-1/4};\zeta),~ (a_1, a_3;\zeta)(a_2,\zeta;\zeta)$\\
\hline $\Phi_4(222)d_2$ & $a_1,a_2,a_3$ & $\zeta\in k$ &
$(a_1,\zeta;\zeta),~ (\zeta^{-1} a_2, a_3;\zeta),~ (a_1, a_3;\zeta)(a_2,\zeta;\zeta)$\\
\hline $\Phi_4(222)e_0$ & $a_1,a_2,a_3$ & $\zeta\in k$ &
$(a_3,\zeta;\zeta),~ (a_2, a_3;\zeta)(a_1,\zeta;\zeta),~ (a_1, a_3;\zeta)(a_2,\zeta^\nu;\zeta)$\\
\hline $\Phi_4(222)e_r$ & $a_1,a_2,a_3$ & $\zeta\in k$ &
$(a_3,\zeta;\zeta),~ (a_2, \zeta a_3;\zeta)(a_1,\zeta^\kappa;\zeta),~ (a_1, a_3;\zeta)(a_2,\zeta;\zeta)$\\
\hline $\Phi_4(2211)g$ & $a_1,a_2,a_3,a_4$ & $\zeta\in k$ &
$(a_2, a_3;\zeta)(a_4,\zeta;\zeta),~ (\zeta^{-1} a_1, a_3;\zeta)$\\
\hline $\Phi_4(2211)h$ & $a_1,a_2,a_3,a_4$ & $\zeta\in k$ &
$(a_2, a_3;\zeta)(a_4,\zeta;\zeta),~ (a_1,\zeta a_3;\zeta)$\\
\hline $\Phi_4(2211)i$ & $a_1,a_2,a_3,a_4$ & $\zeta\in k$ &
$(a_2, a_3;\zeta)(a_4,\zeta;\zeta),~ (a_1,a_3;\zeta)(a_2, \zeta;\zeta)$\\
\hline $\Phi_4(2211)j_1$ & $a_1,a_2,a_3$ & $\zeta\in k$ &
$(a_3,\zeta;\zeta),~ (a_2, a_3;\zeta),~ (a_1,\zeta a_3;\zeta)$\\
\hline $\Phi_4(2211)j_2$ & $a_1,a_2,a_3$ & $\zeta\in k$ &
$(a_1,\zeta;\zeta),~ (a_2, a_3;\zeta),~ (\zeta^{-1} a_1, a_3;\zeta)$\\
\hline $\Phi_4(2211)k$ & $a_1,a_2,a_3$ & $\zeta\in k$ &
$(a_1,\zeta;\zeta),~ (\zeta^{-1} a_2, a_3;\zeta),~ (a_1, a_3;\zeta)$\\
\hline $\Phi_4(2211)l$ & $a_1,a_2,a_3$ & $\zeta\in k$ &
$(a_3,\zeta;\zeta),~ (a_2, a_3;\zeta),~ (a_1,a_3;\zeta)(a_2,\zeta;\zeta)$\\
\hline $\Phi_4(2211)m$ & $a_1,a_2,a_3$ & $\zeta\in k$ &
$(a_1,\zeta;\zeta),~ (a_2, \zeta a_3;\zeta),~ (a_1,a_3;\zeta)$\\
\hline $\Phi_4(2211)n$ & $a_1,a_2,a_3$ & $\zeta\in k$ &
$(a_1,\zeta;\zeta),~ (a_2, a_3;\zeta),~ (a_1,a_3;\zeta)(a_2,\zeta;\zeta)$\\
\hline $\Phi_4(21^4)d$ & $a_1,a_2,a_3,a_4$ & $\zeta\in k$ &
$(a_2, a_3;\zeta),~ (a_1, a_3;\zeta)(a_4,\zeta;\zeta)$\\
\hline $\Phi_4(21^4)e$ & $a_1,a_2,a_3$ & $\zeta\in k$ &
$(a_3,\zeta;\zeta),~ (a_2, a_3;\zeta),~ (a_1,a_3;\zeta)$\\
\hline $\Phi_4(21^4)f$ & $a_1,a_2,a_3$ & $\zeta\in k$ &
$(a_1,\zeta;\zeta),~ (a_2, a_3;\zeta),~ (a_1,a_3;\zeta)$\\
\hline $\Phi_{12}(2211)a$ & $a_1,a_2,a_3,a_4$ & $\zeta\in k$ &
$(a_2,a_4;\zeta)(a_3,\zeta;\zeta),~ (a_1,\zeta a_3;\zeta)$\\
\hline $\Phi_{12}(2211)c$ & $a_1,a_2,a_3,a_4$ & $\zeta\in k$ &
$(a_2,\zeta a_4;\zeta)(a_1, \zeta;\zeta),~ (a_1,\zeta a_3;\zeta)$\\
\hline $\Phi_{12}(2211)d$ & $a_1,a_2,a_3,a_4$ & $\zeta\in k$ &
$(a_2, a_4;\zeta)(a_1, \zeta;\zeta),~ (a_1, a_3;\zeta)(a_2, \zeta;\zeta)$\\
\hline $\Phi_{12}(2211)e$ & $a_1,a_2,a_3,a_4$ & $\zeta\in k$ &
$(a_2, a_4;\zeta)(a_1, \zeta;\zeta),~ (a_1,\zeta a_3;\zeta)(a_2, \zeta;\zeta)$\\
\hline $\Phi_{12}(2211)f$ & $a_1,a_2,a_3,a_4$ & $\zeta\in k$ &
$(a_2, a_4;\zeta)(a_3, \zeta;\zeta),~ (a_1,\zeta a_3;\zeta)(a_2, \zeta;\zeta)$\\
\hline $\Phi_{12}(2211)g$ & $a_1,a_2,a_3,a_4$ & $\zeta\in k$ &
$(a_2,\zeta a_4;\zeta)(a_3, \zeta;\zeta),~ (a_1,\zeta a_3;\zeta)$\\
\hline $\Phi_{12}(2211)h$ & $a_1,a_2,a_3,a_4$ & $\zeta\in k$ &
$(a_2,\zeta a_4;\zeta)(a_3, \zeta;\zeta),~ (a_1,\zeta a_3;\zeta)(a_4, \zeta;\zeta)$\\
\hline $\Phi_{12}(2211)i$ & $a_1,a_2,a_3,a_4$ & $\zeta\in k$ &
$(a_2,\zeta a_4;\zeta)(a_3, \zeta;\zeta),~ (a_1,\zeta a_3;\zeta)(a_2, \zeta;\zeta)$\\
\hline $\Phi_{12}(21^4)b$ & $a_1,a_2,a_3,a_4$ & $\zeta\in k$ &
$(a_2, a_4;\zeta)(a_1, \zeta;\zeta),~ (a_1,\zeta a_3;\zeta)$\\
\hline $\Phi_{12}(21^4)c$ & $a_1,a_2,a_3,a_4$ & $\zeta\in k$ &
$(a_2, a_4;\zeta)(a_1, \zeta;\zeta),~ (a_1, a_3;\zeta)$\\
\hline $\Phi_{12}(21^4)d$ & $a_1,a_2,a_3,a_4$ & $\zeta\in k$ &
$(a_2, a_4;\zeta),~ (a_1,\zeta a_3;\zeta)(a_2, \zeta;\zeta)$\\
\hline $\Phi_{12}(21^4)e$ & $a_1,a_2,a_3,a_4$ & $\zeta\in k$ &
$(a_2,\zeta a_4;\zeta)(a_1, \zeta;\zeta),~ (a_1,\zeta a_3;\zeta)(a_2, \zeta;\zeta)$\\
\hline
\end{tabular}
\end{table}

\clearpage

\begin{table}[ht]
\caption{Pullbacks of order $p^6$ with abelian quotients: Part
II}\centering
\begin{tabular}{|l|l|l|l|} \hline
Group & $\vert k^\times/k^{\times p}\vert$ & Root of unity & Trivial elements in \Br(k)   \\
\hline\hline $\Phi_{13}(2211)a$ & $a_1,a_2,a_3,a_4$ & $\zeta\in k$ &
$(a_1,a_3;\zeta)(a_2,\zeta a_4;\zeta),~ (a_1,\zeta a_2;\zeta)$\\
\hline $\Phi_{13}(2211)b$ & $a_1,a_2,a_3,a_4$ & $\zeta\in k$ &
$(\zeta^{-1} a_1,a_3;\zeta)(a_2,a_4;\zeta),~ (a_1,\zeta a_2;\zeta)$\\
\hline $\Phi_{13}(2211)c_r$ & $a_1,a_2,a_3,a_4$ & $\zeta\in k$ &
$(a_1,a_3;\zeta)(a_2,\zeta^r a_4;\zeta),~ (a_1,a_2;\zeta)(a_3, \zeta;\zeta)$\\
\hline $\Phi_{13}(2211)d$ & $a_1,a_2,a_3,a_4$ & $\zeta\in k$ &
$(a_1,\zeta a_3;\zeta)(a_2, a_4;\zeta),~ (a_1,a_2;\zeta)(a_3, \zeta;\zeta)$\\
\hline $\Phi_{13}(2211)e_r$ & $a_1,a_2,a_3,a_4$ & $\zeta\in k$ &
$(a_1,a_3;\zeta)(\zeta^{-r} a_2, a_4;\zeta),~ (a_1,\zeta a_2;\zeta)$\\
\hline $\Phi_{13}(2211)f$ & $a_1,a_2,a_3,a_4$ & $\zeta\in k$ &
$(a_1,a_3;\zeta)(\zeta^{-1} a_2, a_4;\zeta),~ (a_1, a_2;\zeta)(a_3, \zeta;\zeta)$\\
\hline $\Phi_{13}(21^4)a$ & $a_1,a_2,a_3,a_4$ & $\zeta\in k$ &
$(a_1,a_3;\zeta)(a_2, a_4;\zeta),~ (a_1,\zeta a_2;\zeta)$\\
\hline $\Phi_{13}(21^4)b$ & $a_1,a_2,a_3,a_4$ & $\zeta\in k$ &
$(a_1,\zeta a_3;\zeta)(a_2, a_4;\zeta),~ (a_1, a_2;\zeta)$\\
\hline $\Phi_{13}(21^4)c$ & $a_1,a_2,a_3,a_4$ & $\zeta\in k$ &
$(\zeta^{-1} a_1,a_3;\zeta)(a_2, a_4;\zeta),~ (a_1, a_2;\zeta)$\\
\hline $\Phi_{13}(21^4)d$ & $a_1,a_2,a_3,a_4$ & $\zeta\in k$ &
$(a_1,a_3;\zeta)(a_2, a_4;\zeta),~ (a_1, a_2;\zeta)(a_3, \zeta;\zeta)$\\
\hline $\Phi_{13}(1^6)$ & $a_1,a_2,a_3,a_4$ & $\zeta\in k$ &
$(a_1,a_3;\zeta)(a_2, a_4;\zeta),~ (a_1, a_2;\zeta)$\\
\hline $\Phi_{15}(2211)a$ & $a_1,a_2,a_3,a_4$ & $\zeta\in k$ &
$(a_1,a_3;\zeta)(a_2,\zeta a_4^g;\zeta),~ (a_1,\zeta a_2;\zeta)(a_3,a_4;\zeta)$\\
\hline $\Phi_{15}(2211)b_{r,s}$ & $a_1,a_2,a_3,a_4$ & $\zeta\in k$ &
$(a_1,\zeta^r a_3;\zeta)(a_2,\zeta^\kappa a_4^g;\zeta),~ (a_1,\zeta a_2;\zeta)(a_3,a_4;\zeta)$\\
\hline $\Phi_{15}(2211)c$ & $a_1,a_2,a_3,a_4$ & $\zeta\in k$ &
$(a_1, a_3;\zeta)(\zeta^ga_2^g, a_4;\zeta),~ (a_1,\zeta a_2;\zeta)(a_3,a_4;\zeta)$\\
\hline $\Phi_{15}(2211)d_r$ & $a_1,a_2,a_3,a_4$ & $\zeta\in k$ &
$(a_1, a_3;\zeta)(\zeta^{-\kappa} a_2^g, a_4;\zeta),~ (a_1,\zeta a_2;\zeta)(a_3,a_4;\zeta)$\\
\hline $\Phi_{15}(21^4)$ & $a_1,a_2,a_3,a_4$ & $\zeta\in k$ &
$(a_1, a_3;\zeta)(a_2^g, a_4;\zeta),~ (a_1,\zeta a_2;\zeta)(a_3,a_4;\zeta)$\\
\hline $\Phi_{15}(1^6)$ & $a_1,a_2,a_3,a_4$ & $\zeta\in k$ &
$(a_1, a_3;\zeta)(a_2^g, a_4;\zeta),~ (a_1,a_2;\zeta)(a_3,a_4;\zeta)$\\
\hline
\end{tabular}
\end{table}

\bigskip

\begin{table}[ht]
\caption{Groups of order $p^6$ with quotient $(C_{p^2})^2$ obtained
by factoring out $\mu_{p^2}$}\centering
\begin{tabular}{|l|l|l|l|} \hline
Group & $\vert k^\times/k^{\times p}\vert$ & Root of unity & Trivial elements in \Br(k)   \\
\hline\hline $\Phi_{14}(42)$ & $a_1,a_2$ & $\zeta_{p^2}\in k$ &
$(a_1,\zeta_{p^2} a_2;\zeta_{p^2})$\\
\hline $\Phi_{14}(321)$ & $a_1,a_2$ & $\zeta_{p^2}\in k$ &
$(a_1,\zeta a_2;\zeta_{p^2})$\\
\hline $\Phi_{14}(222)$ & $a_1,a_2$ & $\zeta_{p^2}\in k$ &
$(a_1,a_2;\zeta_{p^2})$\\
\hline
\end{tabular}
\end{table}

%--------------------------------------Section 5-------------------------------------------------
\section*{Acknowledgements}

I am grateful to an anonymous referee for the helpful suggestions
that improved dramatically the present paper and in particular
Theorem \ref{c3.6}.

\end{document}